\documentclass[10pt]{paper}

\usepackage{amssymb,amsmath,enumerate,theorem,epsfig}
\usepackage{amsfonts}
\usepackage{a4,latexsym,parskip}
\usepackage{hyperref}
\hypersetup{pdfauthor=Matthias Sch"utt}
\hypersetup{pdftitle=max singular fibre}

\newcommand{\Phoch}[1]{\mathbb{P}^{#1}}

\newcommand{\C} [1][]{\mathbb{C}^{#1}}
\newcommand{\Q} [1] []{\mathbb{Q}_{#1}}

\newtheorem{Spezial-Theorem}{Theorem}[section]
\newtheorem{Spezial-Proposition}{Proposition}[section]
\theoremstyle{break} \newtheorem{Theorem}{Theorem}[section]
\newtheorem{Proposition}[Theorem]{Proposition}

\newtheorem{Corollary}[Theorem]{Corollary}

\begin{document}
\setlength{\unitlength}{1cm}

\title{The maximal singular fibres of elliptic K3 surfaces}
\author{Matthias Sch\"utt}
\maketitle

\begin{small}

\textbf{Abstract.}\ We prove that the maximal singular fibres of
elliptic K3 surfaces have type $I_{19}$ and $I_{14}^*$ unless the
characteristic of the ground field is 2. In characteristic 2, the
maximal singular fibres are $I_{18}$ and $I_{13}^*$. The paper
supplements work of Shioda in \cite{Sh03} and \cite{ShDS}.



\end{small}
\vspace{0.3cm}

\section{Introduction}

This note investigates the maximal singular fibres of elliptic K3 surfaces
\[
\pi:X\to\Phoch{1}
\]
in arbitrary characteristic. Here, maximality should be understood
in terms of the number of components. Throughout we will assume
the fibration to have a section (since otherwise we can consider
its Jacobian).

In characteristic 0, the maximal fibres are known to have type $I_{19}$ and $I_{14}^*$ in Kodaira's notation. This follows from the Lefschetz bound
\begin{equation}\label{eq:0}
\rho(X) \leq h^{1,1}(X)=20
\end{equation}
by way of the Shioda-Tate formula \cite[Cor. 5.3]{Sh90}. Here,
$\rho(X)=$ rk $NS(X)$ denotes the Picard number, i.e. the rank of
the N\'eron-Severi group. The existence of such K3 fibrations
follows from the classifications of \cite{MP} and \cite{N}.

In positive characteristic $p$, however, we only have the trivial weaker bound
\begin{equation}\label{eq:b_2}
\rho(X) \leq b_2(X)=22
\end{equation}
involving the second Betti number $b_2(X)$. Note that here
$\rho(X)=21$ is impossible by the work of M. Artin \cite{Artin}.
Hence, the existence of a larger singular fibre (than $I_{19}$ or
$I_{14}^*$) of the elliptic K3 surface $\pi:X\to\Phoch{1}$ already
implies
\[
\rho(X)=22,
\]
again due to the Shioda-Tate formula. In other words, $X$ is a supersingular K3 surface.

This property enabled Shioda in \cite[Rem. (4.3)]{Sh03} to exclude
the fibres of type $I_{21}$ and $I_{16}^*$ which are a priori
maximal possible by equation (\ref{eq:b_2}). In \cite{ShDS}, he
reproved existence and maximality of $I_{19}$ and $I_{14}^*$ over
$\C$ by using Davenport-Stothers triples. The remark in \cite[\S
3]{ShDS} states that the arguments for the maximality stay valid
in characteristic $p$ if $p$ is sufficiently large. In this note,
we will employ a different approach to extend these results to
arbitrary characteristic.

\begin{Theorem}\label{Thm:all}
In every characteristic $p\neq 2$, the maximal singular fibres of elliptic K3 fibrations are of type $I_{19}$ and $I_{14}^*$. In characteristic $2$, they are $I_{18}$ and $I_{13}^*$.
\end{Theorem}

The main claim of the theorem will follow from Theorems
\ref{Thm:20} and \ref{Thm:15*}. For odd characteristic, the proofs
rely on a strong property of supersingular K3 surfaces, encoded in
the Artin invariant. We then only use elementary congruences
involving the height pairing.

In characteristic $2$, we will perform explicit calculations involving the Weierstrass equation. In the additive case, we will also make extensive use of Tate's algorithm to determine the explicit type of the special fibre. It then turns out that the methods already apply to the subsequent fibre types. The corresponding results can be found in Propositions \ref{Prop:19} and \ref{Prop:14*}.

The paper is organized as follows: We first consider
multiplicative fibres and give a complete proof for these. This is
followed by an application to the reduction of the [1,1,1,1,1,19]
fibration (Cor. \ref{Cor:bad}), i.e. the elliptic K3 surface with
singular fibres of types $I_1,\hdots,I_{19}$. The note is
concluded by the proofs for additive fibres, but these will only
be sketched roughly.

\section{Multiplicative case}

The following three sections are devoted to the proof that in arbitrary characteristic there is no elliptic K3 surface
\[
\pi:X\to\Phoch{1}
\]
with a fibre of type $I_{20}$. To achieve this, we will assume on
the contrary that there is such a fibration, and establish a
contradiction. In particular, we require that the characteristic
$p$ is positive due to the Lefschetz bound (\ref{eq:0}) in zero
characteristic.

Since $\rho(X)\neq 21$, the existence of a fibre of type $I_{20}$
(or $I_{21}$) already implies
\[
\rho(X)=22
\]
by the Shioda-Tate formula. Hence $X$ is a supersingular K3 surface. 
Then, we have the following result concerning the N\'eron-Severi lattice $NS(X)$ which goes back to M. Artin \cite{Artin} (for $p=2$, this is due to Rudakov and \v Safarevi\v c \cite{RS}):
\begin{Theorem}[Artin, Rudakov-\v Safarevi\v c]\label{Thm:Artin-inv}
Let $X$ be a supersingular K3 surface over a field of
characteristic $p$. Then
\[\text{discr } NS(X)=-p^{2\sigma_0} \;\;\;\text{ for some }\;\sigma_0\in\{1,\hdots,10\}.\]
\end{Theorem}
Here, $\sigma_0$ is called  the \emph{Artin invariant}.

As a direct application, we can rule out the Mordell-Weil group
$MW(\pi)$ of the K3 surface $\pi:X\to\Phoch{1}$ to be finite, if
there is a fibre of type $I_{20}$. This was pointed out by Shioda
in \cite[Rem. (4.3)]{Sh03}. Otherwise there would be exactly one
further reducible fibre. Since both possible types $I_2$ and $III$
correspond to the root lattice $A_1$ (as $I_{20}$ corresponds to
$A_{19}$), this would give
\begin{equation}\label{eq:finite-MW}
\text{discr } NS(X)=-\frac{(\text{discr } A_1) (\text{discr }
A_{19})}{|MW(\pi)|^2}=-10 \;\;\text{ or }\; -40.
\end{equation}
This is clearly impossible. The same argument directly excludes
fibres of type $I_{21}$. Hence a fibre of type $I_{20}$ is a
priori the maximal possible for an elliptic K3 surface in positive
characteristic.

\begin{Theorem}\label{Thm:20}
There is no elliptic K3 fibration with a singular fibre of type $I_{20}$.
\end{Theorem}

The proof of this theorem will require us to distinguish between
even and odd characteristic. By our above considerations, an
elliptic K3 fibration $\pi:X\to\Phoch{1}$ with an $I_{20}$ fibre
would neccessarily have Mordell-Weil group of rank 1. Hence, the
discriminant of $NS(X)$ could be expressed in terms of the height
of the generator $P$ of $MW(\pi)$ (up to torsion). In odd
characteristic, it will then be an easy exercise in congruences
using the height pairing to derive a contradiction to Theorem
\ref{Thm:Artin-inv}.

For characteristic 2, on the contrary, we will have to perform
explicit calculations involving the Weierstrass equation to verify
Theorem \ref{Thm:20}. It will then turn out that the same
arguments apply to fibres of type $I_{19}$.

\begin{Proposition}\label{Prop:19}
There is no elliptic K3 fibration over a field of characteristic 2 with a singular fibre of type $I_{19}$.
\end{Proposition}

This can be compared to all other characteristics where there is
indeed such a fibration. This comes from the
[1,1,1,1,1,19]-fibration over $\Q$ (unique up to isomorphism)
which was given by Shioda in \cite{Sh03}. By way of reducing mod
$p$, his equation gives rise to fibrations in question for all
characteristics $p>3$. Furthermore, a $\Q(\sqrt{-3})$-isomorphic
model with good reduction at $p=3$ was found in \cite{ST}. In
particular, this shows that the maximal multiplicative fibre given
in Theorem \ref{Thm:all} does occur in all characteristics
different from 2. On the other hand, Proposition \ref{Prop:19}
implies the following

\begin{Corollary}\label{Cor:bad}
The [1,1,1,1,1,19] fibration, considered over any number field, cannot have good reduction at a prime above 2.
\end{Corollary}

This answers the question of \cite{ST} which in the first instance motivated this note.

\section{Proof of Theorem \ref{Thm:20} in odd characteristic}

In this section we are going to give a proof of Theorem \ref{Thm:20} for odd characteristic $p>2$. For this purpose, we assume that there is an elliptic K3 fibration $\pi:X\to\Phoch{1}$ with an $I_{20}$ fibre and establish a contradiction. In this setting, we have already excluded finite $MW(\pi)$ by discussing equation (\ref{eq:finite-MW}). In other words, $NS(X)$ is generated by the components of the $I_{20}$ fibre (which is the only reducible fibre), the 0-section $O$ and the section $P$ which generates the rank one Mordell-Weil group.

More precisely, this generation is only up to the torsion in $MW(\pi)$. However, there can only be $p$-torsion. This is because any torsion section of order $q$ coprime to $p$ gives rise to another elliptic K3 surface by way of the quotient by its translation. Here, this is a priori impossible for $q\neq 5$ due to the resulting fibre types which contradict the Euler number.

To exclude 5-torsion as well, we further need that the moduli problem for $\Gamma_1(5)$ is representable. Hence $\pi$ would factor through the corresponding modular surface $Y_1(5)$. In characteristic $p\neq 5$, this has the configuration of singular fibres [1,1,5,5] (conf. \cite{L}). Thus, this fibration cannot be connected to a K3 surface with a fibre of type $I_{20}$ by way of pull-back. This rules out 5-torsion in $MW(\pi)$ for $p\neq 5$.

Let $< , >$ denote the height pairing as introduced by Shioda in
\cite[Sect.~8]{Sh90}. It induces the structure of a
positive-definite lattice on $MW(\pi)/MW(\pi)_{\text{tor}}$. In
our case, we have
\begin{equation}\label{eq:height}
<P,P> = 4+2 \;(P.O)-\frac{i(20-i)}{20}>0.
\end{equation}
Here, $(P.O)$ is the intersection number of the 0-section $O$ and
the section $P$ in $NS(X)$. In particular, this is an integer, but
we will not need any other information about it. Furthermore, $i$
denotes the component $C_i$ of the $I_{20}$ fibre which $P$ meets;
that is, the components $C_j$ are numbered cyclically (up to
orientation) such that $O$ meets $C_0$.

The height pairing uses the projection onto the orthogonal
complement of the trivial sublattice generated by $O$ and the
vertical divisors in $NS(X)\otimes\Q$. Hence we have, up to an
even power of $p$ as explained,
\begin{equation}\label{eq:discr}
\text{discr } NS(X) = -20 <P,P>.
\end{equation}


Combining equations (\ref{eq:height}) and (\ref{eq:discr}), we
obtain
\begin{equation*}
|\text{discr } NS(X)| = 80+40 \;(\phi(P).O)-i(20-i).
\end{equation*}
Comparing with Theorem \ref{Thm:Artin-inv} and reducing mod 8,
this gives
\begin{equation*}
i^2-4i\equiv p^{2\sigma_0} \mod{8}.
\end{equation*}
This congruence is compatible with the indeterminacy in the
torsion of $MW(\pi)$. In particular, $i$ has to be odd. Inserting
leads to the congruence
\begin{equation*}
-3 \equiv 1 \mod{8}.
\end{equation*}
This gives the required contradiction. Thus, we have proved
Theorem \ref{Thm:20} in odd characteristic.

\section[Proof of Theorem \ref{Thm:20} in characteristic $2$]{Proof of Theorem \ref{Thm:20} in characteristic $\mathbf{2}$}
\label{s:proof-2}

The above argument breaks down in characteristic $2$. Instead, we
will perform an explicit calculation involving the general
Weierstrass equation and the corresponding discriminant to prove
Theorem \ref{Thm:20} in this characteristic.

We start with the general Weierstrass equation
\begin{equation}\label{eq:Weier}
y^2+a_1xy+a_3y=x^3+a_2x^2+a_4x+a_6
\end{equation}
where the $a_i$ are polynomials in $\overline{\mathbb{F}}_2[t]$ of
degree (at most) $2i$. In order to define an honest K3 surface
(instead of a rational surface, the product of two elliptic curves
or a singular surface), we further have to impose some conditions
on the $a_i$ (e.g. there is an $i$ such that deg $a_i>i$), but we
will not go into detail with this here.

In characteristic $2$, the discriminant is given by
\begin{equation}\label{eq:Delta}
\Delta=a_1^4(a_1^2a_6+a_1a_3a_4+a_2a_3^2+a_4^2)+a_1^3a_3^3+a_3^4.
\end{equation}
We emphasize that next to the vanishing order of $\Delta$, we also
have to take the possibility of wild ramification into account.
This can only occur at the additive fibres.

Our next aim is to normalize the Weierstrass equation
(\ref{eq:Weier}). For this purpose we adopt the techniques of
\cite{L}. Since our analysis will depend on the quadratic
polynomial $a_1$, we have to distinguish between three cases:
\begin{enumerate}[(i)]
\item $a_1\equiv 0$,\label{(i)} \item $a_1\not\equiv 0$ is a
square in $\overline{\mathbb{F}}_2[t]$,\label{(ii)} \item $a_1$
has two distinct zeroes (possibly including
$\infty$).\label{(iii)}
\end{enumerate}

In case (\ref{(i)}), all the singular fibres lie above the zeroes
of $a_3$ and therefore are additive. Thus there cannot be any
singular fibre of type $I_n, n>0,$ at all.

In case (\ref{(ii)}), we normalize $a_1$ to become $t^2$. Then,
using successive translations
\[
x\mapsto x+\alpha,\;\; y\mapsto y+\beta\;\; \text{and} \;\;
y\mapsto y+\gamma x,
\]
we can assume that \[ a_3=at+b,\;\; a_4=ct+d\;\;\text{and}\;\;
a_2=t\,\tilde{a}_2,\]
respectively, where $a,b,c,d$ are constants
and $\tilde a_2$ has degree 3. Thus equation (\ref{eq:Weier})
becomes
\begin{equation}\label{eq:t^2}
y^2+t^2xy+(at+b)y=x^3+t\,\tilde a_2 x^2+(ct+d)x+a_6
\end{equation}
with discriminant
\begin{equation}
\begin{matrix}\Delta & = & t^8(t^4a_6+t^2(at+b)(ct+d)+t\,\tilde a_2 (at+b)^2+(ct+d)^2)\\
&& \;\;\; +t^6(at+b)^3+(at+b)^4.\end{matrix}
\end{equation}

Note that we still have enough freedom to move a fibre above
$t_0\neq 0$ to $\infty$ without changing the general shape of
equation (\ref{eq:t^2}). This comes from the translation $s\mapsto
s+\frac{1}{t_0}$ in the hidden parameter $s=\frac{1}{t}$ at
$\infty$.

Hence, we only have to decide whether it is possible for an
$I_{20}$ fibre in equation (\ref{eq:t^2}) to sit above 0 or
$\infty$. We shall see that the same argument works for fibres of
type $I_{19}$. In fact, such a fibre above $\infty$ is a priori
impossible, since $\Delta$ does not contain any term corresponding
to $t^5$.

At first, let us assume that the vanishing order of $\Delta$ at $\infty$ is at least 20: $v_\infty(\Delta)\geq 20$. That is, all terms in $\Delta$ of order greater than 4 in $t$ vanish. Writing $\Delta=\sum d_it^i$, the first non-trivial equations read
\[
d_6=b^3=0,\;\;\; d_7=ab^2,\;\;\; d_8=d^2+a^2b=0,\;\;\;d_9=a^3+b^2 \tilde a_2(0)=0,
\]
so $b=d=a=0$. The further elimination process kills every
coefficient except for $\tilde a_2$. But then the resulting
fibration becomes singular in codimension 1 ($\Delta\equiv 0$).
This rules out $v_\infty(\Delta)\geq 20$.

Similarly, if $v_0(\Delta)>4$, then $d_0=b^4=0$ and $d_4 =a^4=0$.
Hence the reduction at $t=0$ becomes additive. This rules out the
existence of an $I_n$ fibre above 0 for $n>4$. For a brief
description of the resulting additive fibre, we refer to the
discussion in Section \ref{ss:even*}.

The remaining case which we have to consider is (\ref{(iii)}).
Here we normalize $a_1=t$. In this setting, the changes of
variables $x\mapsto x+\alpha$ and  $y\mapsto y+\beta$ lead to the
Weierstrass equation
\begin{equation}\label{eq:t}
y^2+txy+(at^6+b)y=x^3+a_2 x^2+(ct^8+d)x+a_6
\end{equation}
with discriminant
\begin{equation}\label{eq:discr-t}
\begin{matrix}\Delta & = & t^4(t^2a_6+t(at^6+b)(c^8t+d)+a_2 (at^6+b)^2+(ct^8+d)^2)\\
&& \;\;\; +t^3(at^6+b)^3+(at^6+b)^4.\end{matrix}
\end{equation}

Note that the Weierstrass equation (\ref{eq:t}) is symmetric in 0
and $\infty$ in the following sense: It takes the same general
shape for the local parameters $t$ and $s=\frac{1}{t}$. On the
other hand, we cannot move a fibre to 0 or $\infty$ anymore while
preserving the type of equation as given in (\ref{eq:t}). Hence we
have to consider the problem of an $I_{20}$ fibre above 0 (or
equivalently $\infty$) and above $t_0\neq 0,\infty$. After
rescaling, we will assume $t_0=1$.

At first, we shall assume that there is a singular fibre above 0
for the Weierstrass equation (\ref{eq:t}). Here $v_0(\Delta)>0$ is
equivalent to $b=0$, but then the fibre has additive reduction.
Hence multiplicative singular fibres above 0 are not compatible
with equation (\ref{eq:t}).

On the other hand, let us investigate an $I_{20}$ fibre above the fixed point $t_0=1$. In more generality, we set
\[
\Delta=(t+1)^{20} g
\]
with a non-zero polynomial $g=\sum g_it^i$ of degree at most 4,
i.e. $v_1(\Delta)\geq 20$. Here,
\[
(t+1)^{20}\equiv t^{20}+t^{16}+t^4+1\mod{2},
\]
such that already $g_1=g_2=0$ by the absence of $t$ and $t^2$
terms in $\Delta$. In total, an extensive comparison of
coefficients again produces only zeroes except for $a_2$. But
then, the fibration once more becomes singular. This rules out an
$I_{20}$ fibre in case (\ref{(iii)}) and thereby concludes the
proof of Theorem \ref{Thm:20}.

\section[Application  to $I_{19}$ fibres]{Application  to $\mathbf{I_{19}}$ fibres}

The methods used in the last section can directly be applied to
the question of the existence of an elliptic K3 surface with an
$I_{19}$ fibre in characteristic 2. We will prove that there is no
such (Prop.~\ref{Prop:19}).

In the last section, we indicated in all cases but the last that the vanishing order $v(\Delta)\geq 19$ also results either in additive reduction or in a singular fibration. Hence, the same arguments also rule out fibres of type $I_{19}$.

Concerning the last situation, the same explicit elimination as
above can be applied to the more general case of
$v_{t_0}(\Delta)\geq 19$. Since this again gives rise to a
singular fibration, we have thus also completed the proof of
Proposition \ref{Prop:19}.

Together with Theorem \ref{Thm:20}, this concludes the
investigation of the maximal multiplicative fibre for an elliptic
K3 surface in positive characteristic $p$. For $p>2$, this is
$I_{19}$, appearing in the reduction of a suitable $\Q$-model of
the [1,1,1,1,1,19] fibration as indicated after Proposition
\ref{Prop:19}. For $p=2$, the maximal multiplicative fibre has
type $I_{18}$. For instance, this occurs in the purely inseparable
base change of degree 2 of the (extremal) rational elliptic
surface with fibres [1,1,1,9]. In other words, the maximal
multiplicative fibres which Theorem \ref{Thm:all} lists, do exist
in some elliptic K3 fibration.


We shall now return to the [1,1,1,1,1,19] fibration in
characteristic 0, as analyzed in \cite{ST} and \cite{Sh03}. Since
it has a model over $\Q$, we can consider this surface over any
number field $K$. Then Proposition \ref{Prop:19} implies that it
has bad reduction at all the primes of $K$ above 2 (Cor.
\ref{Cor:bad}).

More precisely, we will deduce that at any prime of $K$ above 2
the fibration cannot have reduction with only isolated
ADE-singularities. In order to distinguish from good reduction,
one might compare two cases:

On the one hand, consider the $\Q$-model of the surface and its
reduction mod 19. This has an additional $A_1$-singularity, so
strictly speaking, the reduction at 19 is bad. The five $I_1$
fibres degenerate to two fibres of type $II$ and $III$. Meanwhile,
the $I_{19}$ fibre is preserved. This is impossible in residue
characteristic 2.

On the other hand, we can also consider a possible degeneration of
the $I_{19}$ fibre itself. By way of reducing mod 2 and resolving
singularities, it might become additive. A comparable case
consists of the [1,1,1,1,4,18] fibration over $\Q$ (base change of
degree 2 of [1,1,2,8]): Its reduction mod 2 has only 2 singular
fibres where the $I_{18}$ fibre is preserved and the others
degenerate to one singular fibre of type $I_1^*$.

We will see that such a degeneration of an $I_{19}$ fibre in an
elliptic K3 fibration is impossible in any characteristic. Since
the degenerate fibre has to contain a chain of 18 rational
(-2)-curves, its type can only be $I_n^*$ for some $n>14$. Hence,
the claim will follow from Theorem \ref{Thm:15*}.

\section{Additive case}

This section investigates the maximal additive fibre of elliptic
K3 surfaces. In characteristic 0 or $p>3$, we have an obvious
model with an $I_{14}^*$ fibre. This comes from the $\Q$-model of
the [1,1,1,1,14*] fibration given by Shioda in \cite{Sh03} by way
of good reduction. For $p=3$, we can produce a model with good
reduction at $p$ by twisting over $\Q(\sqrt{3})$ and then
translating $x\mapsto x-\frac{s-2s^3}{3}$. The resulting
Weierstrass equation is
\[
y^2=x^3-s\,(1+2s^2)\,x^2-2s^6(1+s^2)\,x-s^{11}.
\]
Here, we can already see that the special fibre above 0 has type
at least $I_9^*$. A closer exhibition of Tate's algorithm gives
the announced $I_{14}^*$ fibre which has all components defined
over $\Q$.

In characteristic 0, this fibre is clearly maximal. In this
section, we want to reprove this for odd characteristics and even
give a stronger statement for characteristic 2:

\begin{Theorem}\label{Thm:15*}
There is no elliptic K3 surface with a fibre of type $I_n^*$ for $n> 14$.
\end{Theorem}

\begin{Proposition}\label{Prop:14*}
In characteristic 2, the maximal additive fibre for an elliptic K3 surface is $I_{13}^*$.
\end{Proposition}

We shall first prove the theorem in odd characteristic and for
special cases. Again, this will rely on the Artin invariant for
supersingular K3 surfaces. For the subsequent proof in
characteristic 2, we will use Tate's algorithm to determine the
type of the special fibre, employing the approach of Section
\ref{s:proof-2}.

Assume the surface had a fibre of type $I_{15}^*$ (or $I_{16}^*$
which is a priori maximal). The Shioda-Tate formula predicts the
supersingularity of the K3 surface as before. Hence we can use the
Artin invariant in order to establish a contradiction. This works
for odd characteristic and two special cases.

For instance, if there was a fibre of type $I_{16}^*$, then
\[
\text{discr } NS(X) = -\text{discr } D_{20} = -4.
\]
Hence, the surface could only live in characteristic 2 by virtue
of Theorem \ref{Thm:Artin-inv}. On the other hand, it would be
extremal (i.e. the Picard number $\rho(X)$ is maximal and the
Mordell-Weil group finite). But this contradicts the
classification of extremal elliptic K3 surfaces in characteristic
2 as given in \cite[Table 1]{I}.

Similarly, if there was an $I_{15}^*$ fibre, but finite $MW(\pi)$, then
\[
\text{discr } NS(X) = -\frac{(\text{discr } D_{19})\, (\text{discr
} A_1)}{|MW(\pi)|^2} = -2 \;\; \text{ or } \; -8
\]
gives a direct contradiction. These special cases were pointed out by Shioda in \cite[Rem. (4.3)]{Sh03}.

The remainining case of Theorem \ref{Thm:15*} consists of an
$I_{15}^*$ fibre together with Mordell-Weil group of rank 1. For
the proof, we will distinguish between even (positive) and odd
characteristic.

\subsection{Odd characteristic}

If the characteristic $p$ is odd, we can again use a congruence
argument involving the Artin invariant and the height pairing to
prove Theorem \ref{Thm:15*}. Here we write $P$ for the generator
of $MW(\pi)$ up to torsion. Note that, as before, there can
possibly only be $p$-torsion.

Since all other fibres are irreducible, the height of $P$ is given by
\begin{equation}\label{eq:height*}
<P,P>=4+2\,(P.O)-\begin{cases}\hspace{0.3cm}0\\\hspace{0.3cm}1\\1+\frac{15}{4}\end{cases}
\end{equation}
depending on the component of the $I_{15}^*$ fibre which $P$ meets
(cf.~\cite[Sect.~8]{Sh90}).

Up to an even power of $p$ which comes from the possible
$p$-torsion in $MW(\pi)$, the discriminant of the N\'eron-Severi
lattice is given by
\[
\text{discr } NS(X)=-4<P,P>.
\]
In the first two cases of (\ref{eq:height*}), this is an even
number. This gives the required contradiction to Theorem
\ref{Thm:Artin-inv}. In the third case of (\ref{eq:height*}), we
obtain
\[
|\text{discr } NS(X)|=16+8 \,(P.O)-4-15.
\]
Modulo 8, this becomes
\[
|\text{discr } NS(X)|\equiv -3 \mod{8}
\]
which again contradicts Theorem \ref{Thm:Artin-inv}. As this
argument is compatible with squares of $p$ (possibly coming from
torsion of $MW(\pi)$), it proves Theorem \ref{Thm:15*} in odd
characteristic.

\subsection{Characteristic 2}
\label{ss:even*}

For the proof of Theorem \ref{Thm:15*} in characteristic $p=2$, we return to the general Weierstrass form (\ref{eq:Weier}) and the corresponding discriminant (\ref{eq:Delta}). In order to prove Proposition \ref{Prop:14*} as well, we will deal with $I_{14}^*$  and $I_{15}^*$ fibres simultaneously. Hence, choosing the fibre to sit above 0, we have the same condition
\[
v_0(\Delta)\geq 20
\]
as investigated in the multiplicative case. In view of the
previous calculations, we have to consider the following three
cases:
\begin{enumerate}[(i)]
\item $a_1\equiv 0$\label{(i*)}
\item $a_1=t^2$,\label{(ii*)}
\item $a_1=t$.\label{(iii*)}
\end{enumerate}
To determine the type of the special fibre at $t=0$, we apply Tate's algorithm. Roughly speaking, this performs suitable changes of variables in order to obtain the greatest simultaneous vanishing orders of $a_3, a_4$ and $a_6$ possible at the fixed point. These respective orders then predict the type of the special fibre. We refer the reader to \cite[IV.9]{Si} for the details and the notation employed.

In case (\ref{(i*)}), we have $\Delta=a_3^4$, so $a_3=t^5$ or
$a_3=t^6$ after normalizing. Note that a change of variables
possibly involved in performing Tate's algorithm, does not affect
$a_3$. Hence the algorithm definitely terminates when $a_{3,5}$
resp.~$a_{3,6}$ enters (i.e.~at vanishing order $v_0(a_3)=5$
resp.~6). The corresponding fibre types are $I_7^*$ and $I_9^*$.

In case (\ref{(ii*)}), a fibre of type $I_n^*$ above 0 with $v_0(\Delta)\geq 20$ is easily seen to require the Weierstrass equation
\[
y^2+t^2xy=x^3+t\tilde a_2x^2+t^8\tilde a_6
\]
with deg $\tilde a_2\leq 3$, deg $\tilde a_6\leq 4$. Here, $t\nmid\tilde a_2$, since otherwise the surface would be rational. A careful analysis shows that the special fibre becomes maximal when $\tilde a_6=et^4$. Here, $e\neq 0$, since otherwise the fibration would be singular ($\Delta\equiv 0$). The change of variable $y\mapsto y+\sqrt{e}\,t^6$ gives the equivalent Weierstrass equation
\[
y^2+t^2xy=x^3+t\tilde a_2x^2+t^8\sqrt{e}\,x.
\]
This shows that the algorithm terminates at type $I_{12}^*$ where
$a_{4,8}=\sqrt{e}$ enters.

In case (\ref{(iii*)}), we start with Weierstrass equation
(\ref{eq:t}) and discriminant (\ref{eq:discr-t}). First of all,
the vanishing order $v_0(\Delta)\geq 20$ implies $b=0$. On the
other hand, $a\neq 0$, since otherwise $\tilde a_6=0$ and the
fibration would be singular. Hence we can normalize such that
$a=1$. We obtain $d=0$ and 14 further equations such that the
Weierstrass equation (\ref{eq:t}) finally reads
\[
y^2+txy+t^6y=x^3+(et^4+ct^3+\tilde a_6)x^2+ct^8x+t^{10}\tilde a_6.
\]
The special fibre can be seen to have type $I_{12}^*$ unless
$c=\sqrt{e}$. Otherwise, it is $I_{13}^*$. Note that
$v_0(\Delta)=21$ if and only if $c=e$. This completes the proofs
of Theorem \ref{Thm:15*} and Proposition \ref{Prop:14*}. In
particular, the maximal additive fibre for an elliptic K3
fibration in characteristic 2 as given in Proposition
\ref{Prop:14*} does exist.

\vspace{0.3cm}

As a corollary, we deduce that no reduction of the [1,1,1,1,1,19]
fibration gives an elliptic K3 surface with degenerate $I_{19}$
fibre. This answers the question stated at the end of the previous
section.

Furthermore, Proposition \ref{Prop:14*} implies that the
[1,1,1,1,14*] fibration has bad reduction at 2. This, however, was
clear a priori because the N\'eron-Severi group of this surface
has discriminant -4.

\vspace{0.8cm}

\textbf{Acknowledgement:} I am indebted to K. Hulek for his
continuous interest and advice. I would also like to thank R.
Kloosterman and the referee for various helpful comments. Partial
support by the DFG-Schwerpunkt 1094 "Globale Methoden in der
komplexen Geometrie" is gratefully acknowledged.

\vspace{0.4cm}

Matthias Sch\"utt\\ Institut f\"ur Algebraische Geometrie\\ Fakult\"at f\"ur Mathematik und Physik\\ Universit\"at Hannover\\
Welfengarten 1\\ 30167 Hannover\\ Germany\\
\texttt{schuett@math.uni-hannover.de}

\newpage

\end{document}